\documentclass[12pt]{article}
\usepackage[T1]{fontenc}
\usepackage[latin1]{inputenc}

\usepackage{amsfonts}
\usepackage{amssymb}
\usepackage{amsmath}
\usepackage{eufrak}
\usepackage[dvips]{graphicx}
\usepackage{latexsym}
\usepackage{verbatim}
\newcommand{\ind}{\makebox[1em]{\raisebox{-.5ex}[0ex][0ex]{\makebox[0em]%
{$\smile$}}\raisebox{.4ex}[0ex][0ex]{\makebox[-.02em]{$|$}}}}

\def\Ind#1#2{#1\setbox0=\hbox{$#1x$}\kern\wd0\hbox to 0pt{\hss$#1\mid$\hss}
\lower.9\ht0\hbox to 0pt{\hss$#1\smile$\hss}\kern\wd0}
\def\ind{\mathop{\mathpalette\Ind{}}}
\def\Notind#1#2{#1\setbox0=\hbox{$#1x$}\kern\wd0\hbox to 0pt{\mathchardef
\nn=12854\hss$#1\nn$\kern1.4\wd0\hss}\hbox to
0pt{\hss$#1\mid$\hss}\lower.9\ht0 \hbox to
0pt{\hss$#1\smile$\hss}\kern\wd0}

\def\thind{\mathop{\mathpalette\Ind{}}^{\text{\th}} }

\def\starind{\mathop{\mathpalette\Ind{}}^{*} }

\newcommand {\thorn} {\text{\th}}

\newcommand{\C}{{\EuFrak C}}

\newcommand{\bi}{\begin{itemize}}
\newcommand{\ei}{\end{itemize}}

\newcommand{\tp}{\operatorname{tp}}
\newcommand{\Th}{\operatorname{Th}}

\newcommand{\Char}{\operatorname{char}}
\newcommand{\graph}{\operatorname{graph}}
\newcommand{\acl}{\operatorname{acl}}
\newcommand{\Gal}{\operatorname{Gal}}

\def\uth{\text{U}^{\text{\th}} }

\newtheorem{theorem}{Theorem}[section]
\newtheorem{lemma}[theorem]{Lemma}
\newtheorem{fact}[theorem]{Fact}
\newtheorem{claim}{Claim}

\newtheorem{proposition}[theorem]{Proposition}
\newtheorem{definition}[theorem]{Definition}
\newtheorem{remark}[theorem]{Remark}

\newtheorem{question}[theorem]{Question}

\newtheorem{main theorem}{Theorem}
\newtheorem{main question}[main theorem]{Question}
\newtheorem{main conjecture}[main theorem]{Conjecture}
\newtheorem{main corollary}[main theorem]{Corollary}

\title{Superrosy fields and valuations}
\author{Krzysztof Krupi\'nski\footnote{Research supported by NCN grant 2012/07/B/ST1/03513}}
\date{}

\begin{document}
\maketitle
\begin{abstract} We prove that every non-trivial valuation on an infinite superrosy field of positive characteristic has divisible value group and algebraically closed residue field. In fact, we prove the following more general result. Let $K$ be a field such that for every finite extension $L$ of $K$ and for every natural number $n>0$ the index $[L^*:(L^*)^n]$ is finite and, if $\Char(K)=p>0$ and $f\colon L \to L$ is given by $f(x)=x^p-x$, the index $[L^+:f[L]]$ is also finite. Then either there is a non-trivial definable valuation on $K$, or every non-trivial valuation on $K$ has divisible value group and, if $\Char(K)>0$, it has algebraically closed residue field. In the zero characteristic case, we get some partial results of this kind.

We also notice that minimal fields have the property that every non-trivial valuation has divisible value group and algebraically closed residue field.
\end{abstract}
%%%Krzys: Dodalem ponizsze dane.

\section{Introduction}

A motivation for our work comes from some open structural questions concerning fields in various model-theoretic contexts. 

A fundamental theorem says that each infinite superstable field is algebraically closed \cite{Ma,ChSh}. An important generalization of superstable theories is the class of supersimple theories and yet more general class of superrosy theories. Superrosy theories with NIP (the non independence property) also form a generalization of superstable theories which is ``orthogonal'' to supersimple theories in the sense that each supersimple theory with NIP is superstable. It is known from \cite{Hr} that perfect PAC (pseudo algebraically closed) fields with small absolute Galois group (i.e. with absolute Galois group possessing only finitely many closed subgroups of every finite index) are supersimple. A well-known conjecture predicts the converse:

\begin{main conjecture}\label{con0}
Each infinite supersimple field is perfect PAC with small absolute Galois group. 
\end{main conjecture}

A complementary conjecture on infinite superrosy fields with NIP  was formulated in \cite{EKP}.
\begin{main conjecture}\label{con1}
Each infinite superrosy field with NIP is either algebraically or real closed.
\end{main conjecture}

Recall that both algebraically closed and real closed fields are superrosy with NIP. After dropping the NIP assumption, one has to extend the list of possibilities in the conclusion of the above conjecture. Namely, since perfect PAC fields with small absolute Galois group as well as orderable PRC (pseudo real closed) fields  with small absolute Galois group are known to be superrosy \cite[Appendix A]{On}, the following conjecture is strongest possible. (See Section 4 for the definition of PRC fields, which is chosen so that PAC fields are PRC).

\begin{main conjecture}\label{con2}
Each infinite superrosy field is perfect  PRC with small absolute Galois group.
\end{main conjecture}

It is known that a PAC field is simple if and only if its absolute Galois group is small \cite{Hr, Ch,ChPi}; it is supersimple if and only if it is perfect and has small absolute Galois group.
%From this and from \cite[Apendix A]{On}, it follows easily that a PRC field is rosy if and only if its absolute Galois group is small. 
Similarly, a PRC field is superrosy if and only if it is perfect and its absolute Galois group is small (see Fact \ref{fact -1}).
Thus, in Conjectures \ref{con0} and \ref{con2}, once we know that the field is PAC [resp. PRC], the rest of the conclusion is automatically satisfied. It is also easy to see that Conjecture \ref{con2} implies Conjecture \ref{con0}, because one can show that orderable PRC fields have strict order property, and so they are not simple (see Remark \ref{remark -2}). By Fact \ref{fact -3}, Conjecture \ref{con2} also implies Conjecture \ref{con1}.

There are also interesting questions and conjectures concerning NIP fields (without assuming superrosiness). By \cite{KSW}, infinite NIP fields are closed under Artin-Schreier extensions.  A. Hasson and S. Shelah formulated some dichotomies between nice algebraic properties of the field in question and the existence of non-trivial definable valuations. In particular, one can expect that the following is true.

\begin{main conjecture}\label{con3}
 Suppose $K$ is an infinite field with NIP with the property that for every finite extension $L$ of $K$ and for every natural number $n>0$ the index $[L^*:(L^*)^n]$ is finite. Then either there is a  non-trivial definable valuation on $K$, or $K$ is either algebraically or real closed.
\end{main conjecture}
%
%In fact, one could formulate a stronger version of this conjecture with the phrase `any elementary extension of' removed from the assumption. 
Note that if a pure field $K$ is algebraically or real closed, then there is no non-trivial definable valuation on $K$ (e.g. because $K$ is superrosy and we have Fact \ref{fact rosy 1}). Notice also that by Facts \ref{fact rosy 1} and \ref{fact rosy 2}, Conjecture \ref{con3} implies Conjecture \ref{con1}.

Another interesting problem is to classify strongly dependent fields \cite[Section 5]{Sh}.

Independently of the questions of A. Hasson and S. Shelah in the NIP context, our approach to attack Conjectures \ref{con1} and \ref{con2} was to assume that the field in question does not satisfy the conclusion and try to produce a non-trivial definable valuation (existence of which contradicts rosiness by Fact \ref{fact rosy 1}). 
%Using this approach, we obtained the following result whose assumptions generalize the situations from Conjectures \ref{con0}, \ref{con1}, \ref{con2} and \ref{con3}, but whose conclusion is weaker than the conclusions of these conjectures. So, one could say that  it is a common approximation of these conjectures.
This approach led us to the following conjecture whose assumptions generalize the situations from Conjectures \ref{con0}, \ref{con1}, \ref{con2} and \ref{con3}, but whose conclusion is weaker than the conclusions of these conjectures (see Section 4 for explanations). So, one could say that  it is a common approximation of these conjectures.

\begin{main conjecture}\label{main conjecture}
Let $K$ be a field such that for every finite extension $L$ of $K$ and for every natural number $n>0$ the index $[L^*:(L^*)^n]$ is finite and, if $\Char(K)=p>0$ and $f\colon L \to L$ is given by $f(x)=x^p-x$, the index $[L^+:f[L]]$ is also finite. Then either there is a non-trivial definable valuation on $K$, or every non-trivial valuation on $K$ has divisible value group and either algebraically or real closed residue field.
\end{main conjecture}

Our main result is the proof of Conjecture \ref{main conjecture} in the positive characteristic case. In fact, we will prove the following theorem.

\begin{main theorem}\label{main theorem 1}
Let $K$ be a field such that for every finite extension $L$ of $K$ and for every natural number $n>0$ the index $[L^*:(L^*)^n]$ is finite and, if $\Char(K)=p>0$ and $f\colon L \to L$ is given by $f(x)=x^p-x$, the index $[L^+:f[L]]$ is also finite. Then either there is a non-trivial definable valuation on $K$, or every non-trivial valuation on $K$ has divisible value group and, in the case when $\Char(K)>0$, it has algebraically closed residue field.
\end{main theorem}

By Facts \ref{fact rosy 1} and \ref{fact rosy 2}, one gets the following corollary.

\begin{main corollary}\label{main corollary 1}
%Every non-trivial valuation on a superrosy field has divisible value group and either algebraically or real closed residue field.
Every non-trivial valuation on a superrosy field of positive characteristic has divisible value group and algebraically closed residue field.
\end{main corollary}

Since infinite NIP fields are closed under Artin-Schreier extensions \cite{KSW}, we also get the following corollary.

\begin{main corollary}\label{main corollary 2}
%Suppose $K$ is a field with NIP with the property that for every finite extension $L$ of any elementary extension of $K$ and for every natural number $n>0$ the index $[L^*:(L^*)^n]$ is finite. Then either there is a non-trivial definable valuation on $K$, or every non-trivial valuation on $K$ has divisible value group and either algebraically or real closed residue field.
Suppose $K$ is a field of positive characteristic, satisfying NIP and with the property that for every finite extension $L$ of $K$ and for every natural number $n>0$ the index $[L^*:(L^*)^n]$ is finite. Then either there is a non-trivial definable valuation on $K$, or every non-trivial valuation on $K$ has divisible value group and algebraically closed residue field.
\end{main corollary}

%In the above conjectures and results, the phrase `of any elementary extension of $K$' can be replaced by `of an $\aleph_0$-saturated model of $\Th(K)$'. In fact, an essential part of the proof of Theorem \ref{main theorem 1} works for the stronger version where the phrase `any elementary extension of' is dropped. This will be addressed in the course of the proof. 

The proof of Theorem \ref{main theorem 1} relies on \cite{Ko}, where the appropriate results on the existence of non-trivial definable valuations under the presence of certain multiplicative or additive subgroups were established. In contrast, directly from the definition of minimality, we obtain the following variant of Corollary \ref{main corollary 1} for minimal fields. 

\begin{main theorem}\label{main theorem 2}
Every non-trivial valuation on a minimal field has divisible value group and algebraically closed residue field.
\end{main theorem}
Recall that the famous Podewski's conjecture predicts that each minimal field is algebraically closed \cite{Po}.

%In the context of any result or conjecture discussed above, after adding $\sqrt{-1}$ to the field in question, we do not loose the assumptions, whereas in the conclusions we should drop the real closed [or replace PRC by PAC] field possibility. Then, 
Our results tells us that, in various situations, all non-trivial valuations on the field in question have divisible value groups and algebraically closed residue fields. The ultimate goal is to show that the original field (i.e. the residue field with respect to trivial valuation) is algebraically closed [or PAC]. So a questions arises what information about the field in question can be deduced from the information that for all non-trivial valuations the value groups are divisible and the residue fields are algebraically closed. A discussion and some questions about it are included in the last section.

The paper is constructed as follows. First we give prelimnaries concerning rosy theories and valuations, listing all the facts from \cite{Ko} which are used in the course of the proof of Theorem \ref{main theorem 1}. Section 2 is devoted to the proof of Theorem \ref{main theorem 1} and some partial results concerning Conjecture \ref{main conjecture} in the zero characteristic case. A very short and easy Section 3 is self-contained and yields a proof of Theorem \ref{main theorem 2}. In the last section, we discuss  some facts and questions concerning potential applications of our results to prove the original conjectures.  

%The main part of the work contained in this paper was done in 2008, with some cases in the proof of Theorem \ref{main theorem 1} completed in 2013. 

Katharina Dupont is currently working around Conjecture \ref{con3} in her Ph.D. project under the supervision of Salma Kuhlmann and in collaboration with Assaf Hasson. She has a different approach than the one presented in this paper (although also based on \cite{Ko}). A few details on this are mentioned in the last section.

The author would like to thank Thomas Scanlon for discussions and suggestions concerning superrosy fields during the visit at Berkeley in 2007.

\section{Preliminaries}

\subsection{Valuations}
In this subsection, we list the definitions and facts from \cite{Ko} which will be useful in this paper. But before that, let us recall the definition of valuation and other basic notions. A good reference for fields with valuations is for example \cite{EnPr}. 

\begin{definition}
A valuation on a field $K$ is a surjective map $v \colon K \to \Gamma \cup \{\infty\}$, where $(\Gamma,+)$ is an ordered group and $\Gamma<\infty$, satisfying the following axioms. For all $x,y \in K$:
\begin{enumerate}
\item $v(x)=\infty \Longrightarrow x=0$,
\item $v(xy)=v(x)+v(y)$,
\item $v(x+y) \geq \min \{ v(x),v(y)\}$.
\end{enumerate}
\end{definition}

Let $v$ be a valuation on a field $K$. Define
$${\cal O}_v:= \{ x \in K : v(x) \geq 0\}.$$
This is a valuation ring, i.e., a subring of $K$ such that for any $x \in K$, either $x\in {\cal O}_v$ or $x^{-1} \in {\cal O}_v$. We say that the valuation $v$ is trivial if $\Gamma = \{ 0 \}$; equivalently, ${\cal O}_v = K$. The group ${\cal O}_v^*$ of units of ${\cal O}_v$ equals $\{ x \in K : v(x)=0\}$. Finally,
$${\cal M}_v:=\{ x \in K: v(x)>0\}$$
is a unique maximal ideal of ${\cal O}_v$, and 
$$\overline{K}_v:= {\cal O}_v/{\cal M}_v$$
is called the residue field of $v$.  

With a valuation $v$ we associated its valuation ring ${\cal O}_v$. Conversely, starting from a valuation ring ${\cal O}$, one can define a valuation $v$ so that ${\cal O}_v={\cal O}$, and one can do it in such a way that both operations are inverses of each other (after the appropriate identification of value groups). Namely, having a valuation ring ${\cal O}$, we define $\Gamma:=K^*/{\cal O}^*$ with $x{\cal O}^* +y {\cal O}^*:=xy{\cal O}^*$, and we order it by $x{\cal O}^*\leq y{\cal O}^* \iff y/x \in {\cal O}$; then we define a valuation $v$ by $v(x)=x{\cal O}^* \in \Gamma$.

\begin{definition}
%We say that a valuation $v : K \to \Gamma \cup \{ \infty \}$ on a field $K$ (possibly with an additional structure) is definable if both $\Gamma$ and $\graph(v)$ are interpretable in $K$. By the above discussion, this is  equivalent to the fact that ${\cal O}_v$ is definable in $K$.
We say that a valuation $v : K \to \Gamma \cup \{ \infty \}$ on a field $K$ (possibly with an additional structure) is definable if the ordered group $\Gamma$ is interpretable in $K$ and after the interpretation of $\Gamma$ in $K$, $\graph(v)$ is definable. By the above discussion, this is  equivalent to the fact that ${\cal O}_v$ is definable in $K$.
\end{definition}

If $K \subseteq L$ is a field extension and $w$ is a valuation on $L$, we say that $w$ extends $v$ if $w \!\!\upharpoonright_ K =v$ (after an isomorphic embedding of the value group of $v$ into the value group of $w$); equivalently, ${\cal O}_w \cap K={\cal O}_v$. In such a situation, we write $(K,{\cal O}_v) \subseteq (L,{\cal O}_w)$. Let $\Gamma_v$ be the value group of $v$ and $\Gamma_w$ the value group of $w$. Then $\Gamma_v$ can be treated as a subgroup of $\Gamma_w$, and $e({\cal O}_w/{\cal O}_v):=[\Gamma_w:\Gamma_v]$ is called the ramification index of the extension $(K,{\cal O}_v) \subseteq (L,{\cal O}_w)$. Similarly, since ${\cal M}_w \cap {\cal O}_v={\cal M}_v$, we get that $\overline{K}_v={\cal O}_v/{\cal M}_v  \hookrightarrow {\cal O}_w/{\cal M}_w=\overline{L}_w$, and $f({\cal O}_w/{\cal O}_v):=[\overline{L}_w: \overline{K}_v]$ is called the residue degree of the extension  $(K,{\cal O}_v) \subseteq (L,{\cal O}_w)$.

\begin{fact}\label{fundamental inequality} Whenever $(K,{\cal O}_v) \subseteq (L,{\cal O}_w)$ with $n:=[L:K]$ is finite, one has the following inequality $e({\cal O}_w:{\cal O}_v)f({\cal O}_w:{\cal O}_v) \leq n$.
\end{fact}

When $v$ and $w$ are two valuations on the same field $K$, we say that $w$ is a coarsening of $v$ if ${\cal O}_v \subseteq {\cal O}_w$. In such a situation, ${\cal M}_w \subseteq {\cal M}_v$ is a prime ideal of ${\cal O}_v$. In this way, one gets a 1-1 correspondence between overrings of ${\cal O}_v$ and prime ideals of ${\cal O}_v$.  Going further, one gets a 1-1 correspondence between this set of prime ideals and the set of convex subgroups of $\Gamma_v$ (the value group of $v$):
$$\begin{array}{llll}
\Delta & \mapsto & p_\Delta := \{ x \in K : v(x) > \delta\;\; \mbox{for all}\;\; \delta \in \Delta\}\\
p & \mapsto & \Delta_p := \{ \gamma \in \Gamma : \gamma,-\gamma< v(x)\;\; \mbox{for all}\;\; x \in p\}.
\end{array}
$$
For details on this, see \cite[Chapter 2.3]{EnPr}. Let us add that two valuations $v$ and $w$ on $K$ are said to be comparable if ${\cal O}_v \subseteq {\cal O}_w$ or ${\cal O}_w \subseteq {\cal O}_v$.

From the model-theoretic perspective,  an important question is when there exists a non-trivial definable valuation on a given field $K$. A deep insight into this question is provided in \cite{Ko}. In particular, with an additive [resp. multiplicative] subgroup $T$ the author associates a certain valuation ring, denoted by ${\cal O}_T$, and he gives a complete characterization of when ${\cal O}_T$ is first order definable in $(K,+,\cdot,0,1, T)$. Here, we recall some definitions and results from \cite{Ko} which we will use later.

\begin{definition}
Let $v$ be a valuation on $K$ and $T$ an additive [resp. multiplicative] subgroup.
\begin{enumerate}
\item $v$ is compatible with $T$ if ${\cal M}_v \subseteq T$ [resp. $1+{\cal M}_v \subseteq T$].
\item  $v$ is weakly compatible with $T$ if ${\cal A} \subseteq T$ [resp. $1+{\cal A} \subseteq T$] for some ${\cal O}_v$-ideal ${\cal A}$ with $\sqrt{A} = {\cal M}_v$.
\item $v$ is coarsely compatible with $T$ if it is weakly compatible with $T$ and there is no proper coarsening $w$ of $v$ such that ${\cal O}_w^* \subseteq T$.
\end{enumerate}
\end{definition}

\begin{fact}\cite[Lemma 1.2]{Ko}\label{non weak case}
Let $v$ be a valuation on $K$. If either $T$ is a multiplicative subgroup such that for some $n \in \omega$, $(K^*)^n \subseteq T$ and $(n,\Char(\overline{K}_v))=1$, or $T$ is an additive subgroup, $ \Char(K)=p$ and $\{ x^p-x : x \in K \}\subseteq T$, then $v$ is (fully) compatible with $T$ if and only if $v$ is weakly compatible with $T$. 
\end{fact}

\begin{fact}\cite[Proposition 1.4]{Ko}\label{O_T non-trivial}
For an additive [resp. multiplicative] subgroup $T$ of a field $K$, any two coarsly compatible valuations are comparable, and there is a unique finest coarsly compatible valuation ring of $K$ which we will denote by ${\cal O}_T$. Moreover, ${\cal O}_T$ is non-trivial, whenever $T$ is proper (i.e., $T \ne K$ [resp. $T \ne K^*$]) and admits some non-trivial weakly compatible valuation.
\end{fact}

The author concludes that for any additive [resp. multiplicative] subgroup $T$ of a field $K$ exactly one of the following possibilities holds:

\begin{itemize}
\item {\bf groups case}: {\em there is a valuation $v$ on $K$ such that ${\cal O}_v^* \subseteq T$.}\\
In this case, ${\cal O}_T$ is the only coarsely compatible valuation ring with this property, and all weakly compatible valuations are fully compatible.
\item {\bf weak case}: {\em there is a weakly, but not fully compatible valuation on $K$.}\\
In this case, ${\cal O}_T$ is the only valuation ring with this property; the weakly compatible valuations are the coarsenings of ${\cal O}_T$; there is no valuation $v$ with ${\cal O}_v^* \subseteq T$.
\item {\bf residue case}: {\em all weakly compatible valuations are fully compatible, and there is no valuation $v$ with ${\cal O}_v^* \subseteq T$.}\\
In this case, ${\cal O}_T$ is the finest fully compatible valuation ring.
\end{itemize}

Now, we are going to recall \cite[Theorem 2.5]{Ko}, which will be the main tool in this paper. This is a complete characterization of when, for a given additive [resp. multiplicative] subgroup $T$ of a filed $K$, the ring ${\cal O}_T$ is definable in the language ${\cal L}:=\{+,\cdot,0,1, T\}$. In fact, in our applications we will only use the positive part of this characterization (namely the cases when one has definability).  Before we  formulate the full characterization, we should emphasis that  here by definability in ${\cal L}$ we do not just mean definability in the structure $(K,+,\cdot,0,1,T)$, but the existence of a formula $\varphi(x)$ in ${\cal L}$ which defines ${\cal O}_{T'}$ in every model $(K',+,\cdot, 0,1,T') \equiv (K,+,\cdot,0,1,T)$. 

\begin{fact}\cite[Theorem 2.5]{Ko}\label{Ko main}
Let $K$ be a field with an additive [resp. multiplicative] subgroup $T$. Denote by ${\cal M}_T$ the maximal ideal of ${\cal O}_T$, and by $\overline{T}$ the subgroup induced by $T$ on the residue field of ${\cal O}_T$. Then ${\cal O}_T$ is definable in ${\cal L}:=\{+,\cdot,0,1,T\}$ in the following cases:\\[3mm]
\begin{tabular}{|l|c|c|} \hline
& $T\leq (K,+)$ & $T \leq K^*$ \\ \hline
group case & iff ${\cal O}_T$ is discrete or $(\forall x \in {\cal M}_T) (x^{-1}{\cal O}_T \nsubseteq T)$ & always \\ \hline
weak case & iff ${\cal O}_T$ is discrete & iff ${\cal O}_T$ is discrete \\ \hline
residue case & always & iff $\overline{T}$ is no ordering \\ \hline
\end{tabular}
\end{fact}

\subsection{Rosy theories}

In this paper, we will only need two properties of rosy groups and fields. Although they are a folklore, for the reader's convenience we recall fundamental definitions concerning rosiness and we give proofs of these two properties:

\begin{fact}\label{fact rosy 1}
Let $K$ be a rosy field. Then there is no non-trivial definable valuation on $K$.
\end{fact}

\begin{fact}\label{fact rosy 2}
Let $G$ be a commutative superrosy group. Then for every natural number $n>0$, if $G[n]:=\{ g \in G : g^n=e\}$ is finite, the index  $[G:G^n]$ is also finite, where $G^n$ denotes the subgroup consiting of $n$-th powers. 

In particular, if $K$ is a superrosy field, then for every $n>0$ the index $[K^*:(K^*)^n]$ is finite and, if $char(K)=p \ne 0$, then the image of the function $f\colon K \to K$ defined by $f(x)=x^p-x$ is a subgroup of finite index in $K^+$. Since any finite extension of an elementary extension of $K$ is also superrosy, this holds for all finite extensions of any elementary extension of  $K$, too.
\end{fact}

For details on rosy theories the reader is referred to \cite{Ad, EaOn, On}, and on rosy groups to \cite{EKP}. More information about rosy groups and fields can be found in \cite{Kr1,Kr2}.

In this subsection, we work in $\C^{eq}$ where $\C$ is a monster model of a theory $T$ in a language ${\mathcal L}$.

A motivation to consider rosy theories is the fact that in a sense it is the largest class of theories which allows the application of techniques from stability theory, especially of basic forking calculus. This class contains stable and more generally simple theories as well as o-minimal theories. The definition of rosiness which justifies what we have just said is the following: $T$ is rosy if there is a ternary relation $\starind$ on small subsets of $\C^{eq}$ satisfying all the basic properties of forking independence in simple theories except for the Independence Theorem. Such a relation will be called an independence relation. There is a concrete, particularly useful independence relation in rosy theories, called \th-independence, which we are going to define now.

A formula $\delta (x,a)$ strongly divides over $A$ if it is non-algebraic and the set of formulas  $\{ \delta (x,a') \}_{a'\models
\tp(a/A)}$ is $k$-inconsistent for some $k\in \mathbb{N}$.
We say that $\delta (x,a)$ \th-divides over $A$ if we can
find some tuple $c$ such that $\delta (x,a)$ strongly divides over
$Ac$.
A formula \th-forks over $A$ if it implies a (finite)
disjunction of formulas which \th -divide over $A$.

We say that a type $p(x)$
\th-divides over $A$ if there is a formula implied by  $p(x)$ which
\th-divides over $A$; \th-forking is similarly defined. We say
that $a$ is \th-independent from $b$ over $A$, denoted
$a\thind_A b$, if $\tp\left( a/Ab\right) $ does not \th -fork over
$A$.  

In rosy theories, $\thorn$-independence is the weakest independence relation in the sense that $a \starind_C b$ implies $a \thind_C b$ for any independence relation $\starind$.

By a rosy group [or field] we mean a group [or field], possibly with an additional structure, whose theory is rosy.

Rosy theories can be also defined by means of local \th-ranks.

\begin{definition}
Given a formula $\psi(x)$, a finite set $ \Phi $ of formulas
with object variables $x$ and parameter variables $y$, a finite set of
formulas $\Theta$ in variables $y,z$, and natural number $k>0$, we define the \th$_{\Phi,\Theta, k}$-rank of $\psi$  inductively as follows:
\begin{enumerate}

\item $\thorn_{\Phi,\Theta,k}(\psi) \geq 0$ if $\psi $
is consistent.

\item  For $\lambda $ a limit ordinal, $\thorn_{\Phi,\Theta,k}(\psi)\geq \lambda $ if $\thorn_{\Phi,\Theta,k}(\psi)\geq \alpha $ for all $\alpha
<\lambda $.

\item $\thorn_{\Phi,\Theta,k}(\psi) \geq \alpha +1$ if
there is $ \varphi \in \Phi $, some $\theta(y;z) \in
\Theta$ and parameter $c$ such that

\begin{enumerate}
\item  $\thorn_{\Phi,\Theta,k}(\psi(x)\land \varphi(x;a))\geq \alpha $ for infinitely many $a\models \theta(y;c) $, and

\item  $\left\{ \varphi \left( x;a\right) \right\} _{a \models \theta(y;c)}$ is
$k-$inconsistent.
\end{enumerate}
\end{enumerate}

\noindent Given a (partial) type $\pi(x)$ we define $\thorn_{\Phi,\Theta,k}(\pi(x))$ to be the minimum of $\thorn_{\Phi,\Theta,k}(\psi)$ for $\psi\in\pi(x)$.  %When $\Phi$ and $\Theta$ each contain only one formula, we will write $\thorn_{\varphi,\theta,k}(\psi)$.

\end{definition}

Recall that a theory is rosy if and only if for each $\psi, \Phi, \Theta, k$ as above, the local thorn rank $\thorn_{\Phi,\Theta, k}(\psi)$ is finite. One could prove Fact \ref{fact rosy 1} using this characterization. However, one can give an immediate proof using another characterization of rosiness in terms of the so-called equivalence ranks considered in \cite[Section 5]{EaOn}.

\begin{definition}
Let $\pi(x)$ be a partial type, and let $\Delta$ be a finite set of formulas in variables $x,y,z$.
Define $eq\mbox{-}rk_\Delta(\pi(x))$ as follows:
\begin{enumerate}
\item  $eq\mbox{-}rk_\Delta(\pi(x)) \geq 0$ if $\pi(x)$ is consistent.
 \item For $\lambda$ a limit ordinal, $eq\mbox{-}rk_\Delta(\pi(x))\geq \lambda$ if  $eq\mbox{-}rk_\Delta(\pi(x)) \geq \alpha$ for
all $\alpha < \lambda$.
\item  $eq\mbox{-}rk_\Delta(\pi(x)) \geq \alpha+1$ if there is some equivalence relation $E(x, y)$
defined by $\delta(x, y, c)$ with $\delta(x,y,z) \in \Delta$ and $c\in \C^{eq}$, and there are representatives  $b_i$,  $i < \omega$,
of different equivalence classes, such that $eq\mbox{-}rk_\Delta(\pi(x) \wedge E(x, b_i)) \geq \alpha$.
\end{enumerate}
\end{definition}

From \cite[Section 5]{EaOn}, we know that $T$ is rosy if and only if for every $\Delta$ and $\pi(x)$ as above,  $eq\mbox{-}rk_\Delta(\pi(x))$ is finite.\\

%We say that a valuation $v : K \to \Gamma \cup \{ \infty \}$ on a field $K$ (possibly with an additional structure) is definable if both $\Gamma$ and $\graph(v)$ are interpretable in $K$; equivalently the valuation ring of $v$ is definable in $K$.\\

\noindent
{\em Proof of Fact \ref{fact rosy 1}.}
Suppose for a contradiction that there is a non-trivial definable valuation $v$ on the rosy field $K$. We can assume $K$ is a monster model. For $x,y$ from the sort of $K$ and $z$ from the sort of $\Gamma$ consider the formula 
$$\delta(x,y,z)=(v(x-y) \geq z).$$ 
Then, for $\gamma \in \Gamma$ the formula $\delta(x,y,\gamma)$ defines an equivalence  relation on $K$ which we denote by $E_\gamma(x,y)$. Put $\Delta=\{ \delta(x,y,z)\}$. 

It is enough to show that $eq\mbox{-}rk_\Delta(\delta(x,b,\gamma)) \geq n$ for all $n \in \omega$, $\gamma \in \Gamma$ and $b \in K$, because then $eq\mbox{-}rk_\Delta(\delta(x,b,\gamma))$ are infinite, so $K$ is not rosy.

We argue by induction on $n$. The case $n=0$ is trivial. Suppose $eq\mbox{-}rk_\Delta(\delta(x,b,\gamma)) \geq n$ for all $n \in \omega$, $\gamma \in \Gamma$ and $b \in K$. Consider any $\gamma \in \Gamma$ and $b \in K$. By saturation, there is $\gamma' \in \Gamma$ such that $\gamma' > \gamma$ and $E_{\gamma'}(x,y)$ refines $[b]_{E_\gamma}$ into infinitely many classes, say with representatives $b_i$, $i \in \omega$. Then, taking $c:=\gamma'$ in the definition of $eq\mbox{-}rk_\Delta$, we get $eq\mbox{-}rk_\Delta(\delta(x,b,\gamma)) \geq n+1$.\hfill $\blacksquare$\\

Using $\thind$, we define $\uth$-rank in the same way as U-rank is defined in stable theories by means of $\ind$, namely $\uth$ is a unique function from the collection of all complete types to the ordinals together with $\infty$ with the property that for any
ordinal $\alpha$, $\uth(p)\geq \alpha+1$ if and only if  there is some tuple
$a$ and some type $q \in S(Aa)$ such that $q\supset
p$, $\uth\left( q \right) \geq \alpha$ and $q$
\th-forks over $A$. 
$\uth$-rank in rosy theories has most of the nice properties that $U$-rank has in stable theories, e.g. it satisfies Lascar Inequalities: 
$$\uth(a/b,A)+\uth(b/A)\leq \uth(a,b/A)\leq \uth(a/b,A)\oplus \uth(b/A).$$
Assume $T$ is rosy. If $D$ is an $A$-definable set, then $\uth(D):=\sup\{\uth(d/A): d \in D\}$. Of course, if this supremum is finite, then it is just the maximum. It turns out that if $D$ is a definable group, then the supremum is also attained \cite[Remark 1.20]{EKP}. For $D$ a definable set, $\uth(D)=0$ if and only if $D$ is finite. There are also Lascar inequalities for groups: For definable groups $H\leq G$ we have 
$$\uth(H)+\uth(G/H) \leq \uth(G) \leq \uth(H) \oplus \uth(G/H).$$

We say that $T$ is superrosy if $\uth(p)<\infty$ for every type $p$;
a group [or field] is superrosy if its theory is superrosy.\\

\noindent
{\em Proof of Fact \ref{fact rosy 2}.} Let  $H=G^n$. We claim that $\uth(G)=\uth(H)$. To see this, take $g \in G$ with $\uth(g/\emptyset)=\uth(G)$. Of course, $g^n \in \acl(g)$, so $\uth(g^n/g)=0$. Since $G[n]$ is finite, we get $g \in \acl(g^n)$, and so $\uth(g/g^n)=0$. Thus, by Lascar Inequalities, we get $\uth(g^n)=\uth(g)$, so $\uth(H)=\uth(G)$. Using Lascar Inequalities for groups, we immediately conclude that $\uth(G/H)=0$, so $G/H$ is finite. \hfill $\blacksquare$

%\subsection{NIP}

%In this paper, we only use one consequence of NIP, namely that each field with NIP is closed under Artin-Schreier extensions.

%Let $T$ be a theory. We work in $\C^{eq}$, where $\C \models T$ is a monster model of $T$.

\begin{definition}
We say that a theory $T$ has the NIP if there is no formula $\varphi(x,y)$ and sequence $\langle a_i \rangle_{i < \omega}$ such that for every $w \subseteq \omega$ there is $b_w$ such that $\models \varphi(a_i,b_w)$ iff $i \in w$.
\end{definition}

The main result of \cite{KSW} says that an NIP field $K$ has no Artin-Schreier extensions, i.e., if $\Char(K)=p>0$, then the function $x \mapsto x^p-x$ from $K$ to $K$ is surjective. The proof in general uses some algebraic geometry, but assuming that the image of the function $x \mapsto x^p-x$ is of finite index, this is an immediate consequence of the existence of $(K^+)^{00}$ (i.e., of the smallest type-definable subgroup of $K^+$ of bounded index (wlog we assume here that $K$ is a monster model)). Indeed, one easily checks that $(K^+)^{00}$ is a non-trivial ideal of $K$, so $(K^+)^{00}=K$. Since the image of $x \mapsto x^p-x$ is a definable subgroup of $K^+$ of finite index, we get that it contains $(K^+)^{00}=K$, so it is equal to $K$.

\section{Superrosy fields}

This section is devoted to the proof of Theorem \ref{main theorem 1} and some observations concerning Conjecture \ref{main conjecture} in the zero characteristic case. In fact, we will not use superrosiness, but only the assumption of Theorem \ref{main theorem 1}, which is a consequence of superrosiness by Fact \ref{fact rosy 2}. For a given field $K$ this assumption is the following:

\begin{enumerate}
\item[(A)] For every finite extension $L$ of $K$ and for any natural number $n>0$ the index $[L^*:(L^*)^n]$ is finite and, if $\Char(K)=p>0$ and $f\colon L \to L$ is given by $f(x)=x^p-x$, then the index $[L^+:f[L]]$ is also finite.
\end{enumerate}

%%%%%%%%%%%%%%%%%%%%%%%%%%%%%%%%%%%%%%%%%%%%%%%%%%%%%%%%%%%%%%
\begin{comment}

\begin{enumerate}
\item[(A)] For every finite extension $L$ of any elementary extension of $K$ and for any natural number $n>0$ the index $[L^*:(L^*)^n]$ is finite and, if $\Char(K)=p>0$ and $f\colon L \to L$ is given by $f(x)=x^p-x$, then the index $[L^+:f[L]]$ is also finite.
\end{enumerate}

Suppose $K$ has an elementary extension with a finite extension $L$ contradicting (A). This is witnessed by a certain type over $\emptyset$ in $\Th(K)$. Therefore, every $\aleph_0$-saturated model of $\Th(K)$ has a finite extension which do not satisfy the requirement of (A). So, in Assumption (A), the phrase `of any elementary extension of $K$' can be replaced  by `of an $\aleph_0$-saturated model of $\Th(K)$'.

In fact, an essential part of Theorem \ref{main theorem 1} works under the following weaker version of Assumption (A):
\begin{enumerate}
\item[(A-)] For every finite extension $L$ of $K$ and for any natural number $n>0$ the index $[L^*:(L^*)^n]$ is finite and, if $\Char(K)=p>0$ and $f\colon L \to L$ is given by $f(x)=x^p-x$, then the index $[L^+:f[L]]$ is also finite.
\end{enumerate}
%
We will be precise on which version of this assumption is used in the course of the proof of Theorem \ref{main theorem 1}.

\end{comment}
%%%%%%%%%%%%%%%%%%%%%%%%%%%%%%%%%%%%%%%%%%%%%%%%%%%%%%%%%%%%%%

We start from some preparatory observations.

\begin{remark}\label{rem 2.1}
If $K \subseteq L$ is an algebraic field extension and $v$ is a non-trivial valuation on $L$, then $v\!\! \upharpoonright_K$ is also non-trivial.
\end{remark}
{\em Proof.} Suppose for a contradiction that $v\!\! \upharpoonright_K$ is trivial. Consider any $a \in L$ with $v(a)>0$. Let $P(x) = x^n+a_{n-1}x^{n-1}+ \ldots + a_0 \in K[x]$ be the minimal monic polynomial of $a$ over $K$. Then $v(a^n)>v(a_{n-1}a^{n-1})>\ldots >v(a_0)=0$, so $\infty=v(0)=v(P(a))=v(a_0)=0$, a contradiction. \hfill $\blacksquare$

\begin{remark}\label{rem 2.2} Let $L$ be a finite extension of a field $K$. Suppose that any non-trivial valuation on $L$ has the property that its residue field is algebraically closed. Then the residue field of any non-trivial valuation on $K$ is either algebraically or real closed.
\end{remark}
{\em Proof.} Consider any non-trivial valuation $v$ on $K$. By Chevalley's Extension Theorem, $v$ has an extension to a (non-trivial) valuation $w$ on $L$. By assumption, $\overline{L}_w$ is algebraically closed. Since, by Fact \ref{fundamental inequality}, $\overline{L}_w$ is a finite extension of $\overline{K}_v$, the conclusion follows. \hfill $\blacksquare$

\begin{lemma}\label{lem 2.3}
Let $v$ be a valuation on a field $K$. Let $F$ be a finite extension of the residue field $\overline{K}_v$. Then there is a finite extension $L$ of $K$ and an extension $w$ of $v$ to $L$ such that the field $F$ is isomorphic over $\overline{K}_v$ to the residue field $\overline{L}_w$.
\end{lemma}
{\em Proof.} It is enough to prove it for 1-generated extensions. So, let $F=\overline{K}_v(\overline{\alpha})$ for some $\overline{\alpha} \in F$. Choose  $P(x) \in {\cal O}_v[x]$ a monic polynomial such that $\overline{P}(x)$ is the minimal monic polynomial of $\overline{\alpha}$ over $\overline{K}_v$, where $\overline{P}(x) \in \overline{K}_v[x]$ denotes the polynomial obtained from $P(x)$ by reducing the coefficients modulo ${\cal M}_v$. Consider a root $\beta$ of $P(x)$, put $L:=K(\beta)$ and take any extension of $v$ to a valuation $w$ on $L$. Since the coefficients of $P(x)$ are in ${\cal O}_v$ with the leading coefficient equal to 1, one easily gets that $w(\beta) \geq 0$, so $\beta \in {\cal O}_w$. 

Now, for $\overline{\beta} := \beta + {\cal M}_w \in \overline{L}_w$, we have that $\overline{P}(\overline{\beta})=0$. Since $\overline{P}$ is irreducible, we get 
$$[\overline{L}_w:\overline{K}_v] \geq [\overline{K}_v(\overline{\beta}): \overline{K}_v] = \deg \overline{P} = \deg P \geq [L:K],$$
so we have everywhere equalities by Fact \ref{fundamental inequality}. Therefore, $\overline{L}_w=\overline{K}_v(\overline{\beta})$ which is isomorphic to $F$ over $\overline{K}_v$. \hfill $\blacksquare$\\

The next lemma recalls a standard method of showing that a given field is algebraically closed.

\begin{lemma}\label{lem 2.4}
Let $K$ be a field. If for every finite extension $L$ of $K$ and for every prime number n>0, $(L^*)^n=L^*$ and, in the case when  $\Char(K)=p>0$, the function $f \colon L \to L$ given by $f(x)=x^p-x$ is onto, then $K$ is algebraically closed.
\end{lemma} 
{\em Proof.} By assumption, $K$ is perfect. If it is not algebraically closed, then it has a proper Galois extension $F$ of minimal degree $k>1$. There is an intermediate field $L$ with $\Gal(F/L) \cong {\mathbb Z}_q$ for some prime number $q$. 

If $q=p=\Char (K)$, then Galois theory tells us $F$ is a splitting field over $L$ of a polynomial of the form $x^p-x-a$ for some $a \in L$, but, by assumption, this polynomial has at least one zero in $L$, so all its zeros are in $L$, a contradiction.

Assume $q \ne \Char(K)$. Since $q \le k$ and for a primitive $q$-th root of unity $\zeta_q$ the field $K(\zeta_q)$ is a Galois extension of $K$ of degree less than $q$, by the choice of $k$, we conclude that $\zeta_q \in K$. So $\zeta_q \in L$. By Galois theory, we conclude that $F$ is a splitting field over $L$ of a polynomial of the form $x^q-a$ for some $a \in L$. By assumption, this polynomial has a zero in $L$, so all its zeros are in $L$, a contradiction. \hfill $\blacksquare$

\begin{remark}\label{rem -11}
Let $K$ be an infinite field of characteristic $p>0$. Then either $K^p=K$, or the index $[K^*:(K^*)^p]$ of multiplicative groups is infinite and the index $[K^+:(K^+)^p]$ of additive groups is also infinite. In particular, if the index $[K^*:(K^*)^p]$ is finite, then $K$ is perfect.
\end{remark}
{\em Proof.} Suppose $K^p \ne K$, and let $b_1,b_2 \in K$ be linearly independent over $K^p$. Then $b_1 +k^pb_2$ for $k \in K$ are in different cosets modulo $(K^*)^p$. So the index $[K^*:(K^*)^p]$ is infinite. Similarly, the elements $b_1k^p$ for $k \in K$ are in different additive cosets modulo $(K^+)^p$, so the index $[K^+:(K^+)^p]$ is also infinite \hfill $\blacksquare$\\

\noindent
{\em Proof of Theorem \ref{main theorem 1}.}
Suppose there is no non-trivial definable valuation on $K$.\\[2mm]
{\bf Part 1.} {\em The value group of any non-trivial valuation $v$ on $K$ is divisible.} \\[2mm]
%This is true under (A-).}\\[2mm]
%
{\em Proof.} Suppose for a contradiction that some $\gamma \in \Gamma_v$ is not divisible by $n$. 
Put $${\cal S}:=\{ x(K^*)^n \in K^*/(K^*)^n : {\cal O}_v^* \cap x(K^*)^n \ne \emptyset \}$$ and $$T:=\left\langle \bigcup S \right\rangle,$$ that is, $T$ is the subgroup of $K^*$ generated by the union of the family ${\cal S}$.
\begin{claim} 
$T$ is a proper, definable subgroup of $K^*$ such that ${\cal O}_v^* \leq T$ and $(K^*)^n \leq T$.
\end{claim}
{\em Proof of Claim 1.} Since $(K^*)^n$ is a finite index subgroup of $K^*$, we get that $T$ is a union of finitely many cosets of $(K^*)^n$, so $T$ is definable and contains $(K^*)^n$. The fact that ${\cal O}_v^* \leq T$ follows directly from the definition of $T$. To see that $T$ is proper, suppose for a contradiction that $T=K^*$. Then there is $a \in T$ such that $v(a)=\gamma$. By the definition of $T$, we can write $a=a_1^{\varepsilon_1} \cdot \ldots \cdot a_m^{\varepsilon_m}$, where $a_i \in {\cal O}_v^*\cdot (K^*)^n$ and $\varepsilon_i \in \{-1,1\}$ for $i=1,\dots,m$. Then $\gamma=v(a)= \varepsilon_1 v(a_1)+\ldots + \varepsilon_m v(a_m) \in n\Gamma_v$, a contradiction. \hfill $\blacksquare$\\

By Claim 1, we know that $T$ is proper and $1+{\cal M}_v \subseteq {\cal O}_v^* \leq T$, so $v$ is fully compatible with $T$. By Fact \ref{O_T non-trivial}, we conclude that ${\cal O}_T$ is non-trivial. We also get that we are in the group case, so Fact \ref{Ko main} tells us that ${\cal O}_T$ is definable in $(K,+,\cdot,0,1,T)$, but $T$ is definable in $K$, so ${\cal O}_T$ is definable in $K$. This is a contradiction, which completes the proof of Part 1. \hfill $\blacksquare$\\

\noindent
{\bf Part 2.} %{\em The residue field of any non-trivial valuation on $K$ is either algebraically or real closed.}\\[2mm]
{\em Assuming that $\Char(K)=p>0$, the residue field of any non-trivial valuation $v$ on $K$ is algebraically closed.}\\[2mm]
{\em Proof.} By Lemmas \ref{lem 2.3} and \ref{lem 2.4}, we will be done if we prove the following two claims.

\begin{claim}
For every finite extension $L$ of $K$, for every extension $w$ of $v$ to a valuation on $L$, and for every prime number $n$, one has $(\overline{L}_w^*)^n=\overline{L}_w^*$.
\end{claim}

\begin{claim}
%Assume $\Char(K)=p>0$. Then 
For every finite extension $L$ of $K$, one has that the function $f \colon \overline{L}_w \to \overline{L}_w$ given by $f(x)=x^p-x$ is onto.
\end{claim}

\noindent
{\em Proof of Claim 2.} First, notice that by Remark \ref{rem 2.1}, we can assume that $L=K$. Indeed, since $L$ is definable in $K$ (living in some Cartesian power of $K$), by Remark \ref{rem 2.1}, we get that there is no non-trivial definable valuation on $L$; it is also clear that $L$ satisfies (A).\\[-2mm]

Consider the case $\pmb{n \ne p}$.\\[2mm]
{\bf Subclaim} $K^n(1+{\cal M}_v)=K$\\[2mm]
{\em Proof of the subclaim.}  Suppose it is not true. Then $T:=(K^*)^n(1+{\cal M}_v)$ is a proper subgroup of $K^*$, and $v$ is fully compatible with $T$. By Fact \ref{O_T non-trivial}, ${\cal O}_T$ is non-trivial. 

Since $[K^*:(K^*)^n]$ is finite, we see that $T$ is definable, so there is no non-trivial valuation on $K$ definable in $(K,+,\cdot, 0,1,T)$.

Let $\overline{K}_T$ be the residue field corresponding to the valuation ring ${\cal O}_T$.
Since $(K^*)^n \leq T$ and $(n,\Char(\overline{K}_T))=1$, Fact \ref{non weak case} gives us that we are either in the group case or in the residue case. By Fact \ref{Ko main}, in the group case, ${\cal O}_T$ is definable in $(K,+,\cdot, 0,1,T)$ which is impossible, and, in the residue case,  either ${\cal O}_T$ is definable in $(K,+,\cdot, 0,1,T)$ which is impossible or $\overline{T}$ is the positive cone of an ordering on $\overline{K}_T$. However, the last thing is also impossible as $\Char(\overline{K}_T)=p>0$. \hfill $\blacksquare$\\

By the subclaim, $\overline{K}_v^n=\overline{K}_v$. Indeed, for any $a \in {\cal O}_v$, $a=k^n(1+m)$ for some $k \in K$ and $m \in {\cal M}_v$. Since $1+m \in {\cal O}_v^*$, we get that $k \in {\cal O}_v$, so $\overline{a}=\overline{k}^n \in \overline{K}_v^n$.\\[-2mm]

%Now, consider the case $\pmb{n=p}$. Suppose for a contradiction that $\overline{K}_v^p$ is a proper subfield of $\overline{K}_v$. Let $b_1,b_2 \in \overline{K}_v$ be linearly independent over $\overline{K}_v^p$. Then, the elements $b_1+kb_2$ for $k \in \overline{K}_v^p$ are in different multiplicative cosets modulo $(\overline{K}_v^*)^p$. Hence, $[\overline{K}_v^* : (\overline{K}_v^*)^p]$ is at least the cardinality of $\overline{K}_v^p$, which is infinite, because $\overline{K}_v$ is infinite by the already proved fact that $\overline{K}_v^m=\overline{K}_v$ whenever $(m, \Char(\overline{K}_v))=1$.  Moreover, if $\overline{a}\overline{K}_v^p \ne \overline{b}\overline{K}_v^p$ for some $a,b \in {\cal O}_v$, then $aK^p \ne bK^p$, because otherwise either $ak^p=b$ or $a=bk^p$ for some $k \in {\cal O}_v$, and so either $\overline{a}\overline{k}^p=\overline{b}$ or $\overline{a}=\overline{b}\overline{k}^p$, a contradiction.  This implies that $[K^*:(K^*)^p]$ is infinite, which contradicts (A).\hfill $\blacksquare$\\

Now, consider the case $\pmb{n=p}$. Suppose for a contradiction that $\overline{K}_v^p$ is a proper subfield of $\overline{K}_v$.  As we have already proved that $\overline{K}_v^m=\overline{K}_v$ whenever $(m, \Char(\overline{K}_v))=1$, we see that $\overline{K}_v$ is infinite. Thus, Remark \ref{rem -11} implies that $[\overline{K}_v^* : (\overline{K}_v^*)^p]$ is infinite.  Moreover, if $\overline{a}\overline{K}_v^p \ne \overline{b}\overline{K}_v^p$ for some $a,b \in {\cal O}_v$, then $aK^p \ne bK^p$, because otherwise either $ak^p=b$ or $a=bk^p$ for some $k \in {\cal O}_v$, and so either $\overline{a}\overline{k}^p=\overline{b}$ or $\overline{a}=\overline{b}\overline{k}^p$, a contradiction.  This implies that $[K^*:(K^*)^p]$ is infinite, which contradicts (A).\hfill $\blacksquare$\\

\noindent
{\em Proof of Claim 3.} As in the proof of Claim 2, by Remark \ref{rem 2.1}, we can assume that $K=L$. Let $F \colon K \to K$ be given by $F(x)=x^p-x$.\\[2mm]
{\bf Subclaim} $F[K]+{\cal M}_v=K$.\\[2mm]
{\em Proof of the subclaim.} Suppose it is not true. Then $T:=F[K] + {\cal M}_v$ is a proper subgroup of $K^+$, and $v$ is fully compatible with $T$. By Fact \ref{O_T non-trivial}, ${\cal O}_T$ is non-trivial. As $[K^+:F[K]]$ is finite, $T$ is definable in $K$, so there is no non-trivial valuation on $K$ definable in $(K,+,\cdot, 0,1,T)$.

Since $F[K] \leq T$, Fact \ref{non weak case} ensures that we are either in the group case or in the residue case. In the residue case, ${\cal O}_T$ is definable by Fact \ref{Ko main}, a contradiction. In the group case, since we know that ${\cal O}_T$ is not definable, Fact \ref{Ko main} yields some $x \in {\cal M}_T$ with $x^{-1}{\cal O}_T \subseteq  T$.

Now, we will adopt the argument from the proof of \cite[Theorem 3.1]{Ko}. For the reader's convenience, we give all the details.  

Let ${\cal A}_T$ be the largest fractional ${\cal O}_T$-ideal contained in $T$. We aim at defining a certain fractional ${\cal O}_T$-ideal ${\cal A}_T^\alpha$ which will properly contain ${\cal A}_T$ and which will be contained in $T$. This will contradict the choice of ${\cal A}_T$.

Denote by $v_T$ and by $\Gamma_T$ the valuation and the value group corresponding to ${\cal O}_T$.
For a real number $\alpha$ and $\gamma_1,\gamma_2 \in \Gamma$, by $\gamma_1 \geq \alpha \gamma_2$ we mean that $\gamma_1 \geq r\gamma_2$ for all rationals $r\leq \alpha$ if $\gamma_2\geq 0$, or for all rationals $r \geq \alpha$ if $\gamma_2 <0$. This definition makes sense since $\Gamma_T$ is divisible by Part 1 of the proof. 
By $\gamma_1 < \alpha \gamma_2$ we mean the negation of $\gamma_1 \geq \alpha \gamma_2$. 
For $\alpha>1$ define 
$${\cal A}_T^\alpha:= \{ x \in K : v_T(x) \geq \alpha v_T(y)\;\; \mbox{for some}\;\; y \in {\cal A}_T \}.$$
This is, of course, a fractional ${\cal O}_T$-ideal. Since $x^{-1}{\cal O}_T \subseteq  T$ and $x \in {\cal M}_T$, we have that ${\cal A}_T$ properly contains ${\cal O}_T$. This easily implies that ${\cal A}_T \subseteq {\cal A}_T^\alpha$. Now, we check that this inclusion is proper. 

Take a natural number $n$ such that $\alpha>1+\frac{1}{n}$. There is $y \in {\cal A}_T \setminus {\cal O}_T$ such that $(1+\frac{1}{n})v_T(y) \notin v_T({\cal A}_T)$, as otherwise one easily gets that for all $y \in {\cal A}_T \setminus {\cal O}_T$, $y{\cal O}_T$ generates a valuation ring properly containing ${\cal O}_T$ and contained in $T$, which contradicts the fact that $v_T$ is coarsely compatible with $T$. Choose $z \in K$ with $v_T(z)=\frac{1}{n}v_T(y)$. Then $v_T(zy)=(1+\frac{1}{n})v_T(y) \geq \alpha v_T(y)$, so $zy \in {\cal A}_T^\alpha \setminus {\cal A}_T$.

The proof of the subclaim will be completed if we show that ${\cal A}_T^\alpha \subseteq T$ for $\alpha \in (1,2-\frac{1}{p})$. It is enough to prove that every element $t \in {\cal A}_T^\alpha \setminus {\cal A}_T$ belongs to $T$. We have  $v_T(t)\geq\alpha v_T(y)$ for some $y \in {\cal A}_T\setminus {\cal O}_T$. Since $t \notin {\cal A}_T$, we immediately get $v_T(ty^{-1})=v_T(t)-v_T(y)<0$.
%and, since $\alpha<2$, we have $v_T(ty^{-1}) \geq (\alpha-1)v_T(y)>v_T(y)$. 
%Therefore, $ty^{-1} \in {\cal A}_T \setminus {\cal O}_T$. 
Therefore, $ty^{-1} \in K \setminus {\cal O}_T$. 
Thus, by Claim 2, $ty^{-1} = a^p$ for some $a \in K \setminus {\cal O}_T$, and so $ta^{-p}=y \in {\cal A}_T \setminus {\cal O}_T \subseteq T \setminus {\cal O}_T$, hence $ta^{-p}=b^p-b + m$ for some  $b \in K \setminus {\cal O}_T$ and $m \in {\cal M}_v \subseteq {\cal M}_T$. As $\alpha <2-\frac{1}{p}$, we obtain 
$$\begin{array}{lll}
v_T(a^pb) &= &v_T(a^p)+v_T(b)=v_T(ty^{-1})+v_T(b)= v_T(t)-v_T(y)+ \frac{1}{p}v_T(y) \\ 
& \geq & \left(\alpha -1+\frac{1}{p}\right)v_T(y) > v_T(y).
\end{array}$$
Therefore, $v_T(ab)>v_T(a^pb) > v_T(y)\geq \alpha v_T(y)$ and $v_T(a^pm) >v_T(a^p b) > v_T(y)\geq \alpha v_T(y)$. Hence, $ab, a^pb, a^pm \in {\cal A}_T \subseteq T$, and so $t=((ab)^p -ab)+ ab-a^pb+a^pm \in T$. \hfill $\blacksquare$\\

By the subclaim, $f[\overline{L}_v]=\overline{L}_v$. Indeed, for any $a \in {\cal O}_v$, $a=k^p-k+m$ for some $k \in K$ and $m \in {\cal M}_v$. Then $k \in {\cal O}_v$, so $\overline{a}=\overline{k}^p-\overline{k}$. \hfill $\blacksquare$\\

So, the proof of Part 2 and of the whole theorem has been completed. \hfill $\blacksquare$\\

As was pointed out in the introduction, Corollary \ref{main corollary 1} follows from Theorem \ref{main theorem 1} and Facts \ref{fact rosy 1} and \ref{fact rosy 2}, and Corollary \ref{main corollary 2} follows from Theorem \ref{main theorem 1} and the fact that NIP fields are closed under Artin-Schreier extensions.

\begin{proposition}\label{prop partial}
Let $K$ be a field of characteristic zero satisfying (A) and containing $\sqrt{-1}$. Then either there is a non-trivial definable valuation on $K$, or for every non-trivial valuation $v$ on $K$:
\begin{enumerate}
\item if $\Char(\overline{K}_v)>0$, then $\overline{K}_v^n=\overline{K}_v$ for every $n>0$,
\item if $\Char(\overline{K}_v)=0$, then there is a prime number $p$ such that for all primes $n$ different from $p$, $\overline{K}_v^n=\overline{K}_v$.
\end{enumerate}
\end{proposition}
{\em Proof.} 
Suppose there is no non-trivial definable valuation on $K$.\\[1mm]
(1) The same argument (based on Facts \ref{non weak case}, \ref{O_T non-trivial} and \ref{Ko main}) as in the proof of Claim 2 in the proof of Theorem \ref{main theorem 1} works, noting that the assumption $\sqrt{-1} \in K$ eliminates the possibility that $\overline{T}$ is the positive cone of an ordering on $\overline{K}_T$.\\[1mm]
(2) Suppose that for some prime $p$, $\overline{K}_v^p \ne \overline{K}_v$. Then $T:=(K^*)^p(1 + {\cal M}_v)$ is a proper, definable subgroup of $K^*$. If $p \ne \Char(\overline{K}_T)$, then we get a contradiction as in the proof of the subclaim in the proof of Claim 2 in Theorem \ref{main theorem 1}. So $p=\Char(\overline{K}_T)$. 

Consider any prime $n \ne p$. Let $T':=(K^*)^{np}(1+{\cal M}_v)$. Then $T'\leq T$. On the other hand, since $T$ is definable and ${\cal O}_T$ is not definable, Fact \ref{Ko main} implies that $T$ does not belong to the group case. Thus,  we conclude that ${\cal O}_T \subseteq {\cal O}_{T'}$, and hence $\Char(\overline{K}_{T'}) \in \{ 0, p\}$. However, it is impossible to have $\Char(\overline{K}_{T'})=0$, as in this case once again one gets a contradiction as in the proof of the subclaim in the proof of Claim 2 in Theorem \ref{main theorem 1}. So $\Char(\overline{K}_{T'})=p$. Our goal is to show that $K^n(1 + {\cal M}_v)=K$. If this is not the case, then $T'':=(K^*)^n(1+{\cal M}_v)$ is a proper, definable subgroup of $K^*$, and once again Facts \ref{non weak case}, \ref{O_T non-trivial} and  \ref{Ko main} yield that $(n,\Char(K_{T''}) )\ne 1$ and $T''$ does not belong to the group case. Since $T' \leq T''$, we conclude that ${\cal O}_{T''} \subseteq {\cal O}_{T'}$, but $\Char(\overline{K}_{T'})= p$, and so we get that $\Char(\overline{K}_{T''})=p$. This contradicts the assumption that $n$ is relatively prime to $p$. \hfill $\blacksquare$\\

If one was able to strengthen the conclusion of Proposition \ref{prop partial}(2) by showing that  $\overline{K}_v^n=\overline{K}_v$ for all primes $n$, then using Lemma \ref{lem 2.4}, one would get that, under the assumption of Proposition \ref{prop partial}, either there is a non-trivial definable valuation on $K$, or for any non-trivial valuation $v$ with $\Char(\overline{K}_v)=0$ the residue field $\overline{K}_v$ is algebraically closed. So, arguing as in the proof of Remark \ref{rem 2.2}, one would also get that Assumption (A) implies that either there is a non-trivial definable valuation on $K$, or for any non-trivial valuation $v$ with $\Char(\overline{K}_v)=0$ the residue field $\overline{K}_v$ is either algebraically or real closed. In fact, strengthening slightly Assumption (A), this would imply the full conclusion of Conjecture \ref{main conjecture}.

\begin{proposition}\label{prop hypothetical}
Suppose that the conclusion of Proposition \ref{prop partial}(2) can be strengthened to ` $\!\overline{K}_v^n=\overline{K}_v$ for all primes $n$'.
Let $K$ be a field such that for every finite extension $L$ of any elementary extension of $K$ and for every natural number $n>0$ the index $[L^*:(L^*)^n]$ is finite and, if $\Char(K)=p>0$ and $f\colon L \to L$ is given by $f(x)=x^p-x$, the index $[L^+:f[L]]$ is also finite. Then either there is a non-trivial definable valuation on $K$, or every non-trivial valuation on $K$ has divisible value group and either algebraically or real closed residue field.
\end{proposition}
{\em Proof.} Suppose there is no non-trivial definable valuation on $K$. Using Remarks \ref{rem 2.1} and \ref{rem 2.2}, we can assume that $\sqrt{-1} \in K$. Let $v$ be any non-trivial valuation on $K$. Divisibility of $\Gamma_v$ was proved in Theorem \ref{main theorem 1}. By this theorem, it remains to consider the case $\Char(K)=0$. 
%Using Remarks \ref{rem 2.1} and \ref{rem 2.2}, we can assume that $\sqrt{-1} \in K$.  
From the discussion above Proposition \ref{prop hypothetical}, we are done in the case $\Char(\overline{K}_v)=0$. So it remains to consider the case $\Char(\overline{K}_v)=p>0$. 

Take a monster model $(\widetilde{K},\widetilde{\Gamma},\widetilde{v}) \succ (K,\Gamma_v,v)$.
 Let ${\cal O}$ be a maximal valuation ring in $\widetilde{K}$ containing ${\cal O}_{\widetilde{v}}$ and such that $\frac{1}{p} \notin {\cal O}$. Let ${\cal M}$ be the maximal ideal of ${\cal O}$. Then $p \in {\cal M}$, so $k := {\cal O}/{\cal M}$ is of characteristic $p$.

We claim that there is a non-trivial valuation ring ${\cal O}_1$ in $\widetilde{K}$ such that ${\cal O} \subsetneq {\cal O}_1$. To see this, recall that ${\cal O}$ corresponds to the proper convex subgroup $\widetilde{\Gamma}_{\cal O}$ of $\widetilde{\Gamma}$ defined by
$$\widetilde{\Gamma}_{\cal O}:= \{ \gamma \in \widetilde{\Gamma}: (\forall x \in {\cal M})(\widetilde{v}(x) > \gamma, -\gamma)\}.$$
So, there is $\gamma > \widetilde{\Gamma}_{\cal O}$. Put $\widetilde{\Gamma}_1:= \bigcup_n (-n\gamma,n\gamma)$ a convex subgroup of $\widetilde{\Gamma}$ properly containing $\widetilde{\Gamma}_{\cal O}$. By saturation, $\widetilde{\Gamma}_1 \ne \widetilde{\Gamma}$. So, $\widetilde{\Gamma}_1=\widetilde{\Gamma}_{{\cal O}_1}$ for some non-trivial valuation ring ${\cal O}_1 \supsetneq {\cal O}$.

Since ${\cal O} \subsetneq {\cal O}_1$, we have that $\frac{1}{p} \in {\cal O}_1$, hence the residue field $k_1:={\cal O}_1/{\cal M}_1$ (where ${\cal M}_1$ is the maximal ideal of ${\cal O}_1$) is of characteristic $0$. So $k_1$ is algebraically closed (as the conclusion of Proposition \ref{prop hypothetical} holds in the residue zero characteristic case, also for $\widetilde{K}$ in place of $K$). As ${\cal O}_{\widetilde{v}} \subseteq {\cal O}_1$, this implies that the residue field $\overline{\widetilde{K}}_{\widetilde{v}}$ is algebraically closed, and so $\overline{K}_v$ is algebraically closed, too. \hfill $\blacksquare$

\section{Minimal fields}
We will prove here Theorem \ref{main theorem 2}. In contrast to Theorem \ref{main theorem 1}, where the proof relies on non-trivial results from \cite{Ko}, the proof of Theorem \ref{main theorem 2} is a trivial consequence of the definition of minimality.

Recall that a minimal field is an infinite field whose every definable (in one variable) subset is finite or co-finite. \\

\noindent
{\em Proof of Theorem \ref{main theorem 2}.}  Consider any non-trivial valuation $v$ on $K$.

By minimality, for every natural number $n>0$, $(K^*)^n=K^*$.
%and, if $\Char(K)=p>0$, then the function $f\colon K \to K$ given by $f(x)=x^p-x$ is onto. 
(To see this, notice that $(K^*)^n$ 
%[and $f[K]$ if $\Char(K)=p$] 
is an infinite multiplicative 
%[additive] 
subgroup of $K$, so it is a co-finite subgroup, so it is everything.) From this, it is clear that $\Gamma_v$ is divisible.

Consider any monic polynomial $P(x)=x^n+\overline{a}_{n-1}x^{n-1}+\ldots+\overline{a_0} \in \overline{K}_v$ of positive degree $n$; here $a_0,\dots,a_{n-1} \in {\cal O}_v$. Let $P(x)=x^n+a_{n-1}x^{n-1}+\ldots+a_0$. Since $P(x)$ takes co-finitely many values in $K$ and ${\cal M}_v$ is infinite, there exists $a \in K$ such that $P(a) \in {\cal M}_v$. Then $a \in {\cal O}_v$, as otherwise $v(a)<0$, and so $v(a_ia^i)=v(a_i)+ iv(a) > nv(a)=v(a^n)$ for $i=0,\dots,n-1$, which implies that $v(P(a)) =v(a^n)=nv(a)<0$, a contradiction with the fact that $P(a) \in {\cal M}_v$. Therefore, $\overline{P}(\overline{a}) =0$, i.e., $\overline{P}(x)$ has a root in $\overline{K}_v$. \hfill $\blacksquare$

\section{Final comments}

We will say that a given field $K$ has Property $(*)$ if every non-trivial valuation on $K$ has divisible value group and algebraically closed residue field. We will say that it has Property $(*-)$ if every non-trivial valuation on $K$ has divisible value group and either algebraically or real closed residue field. By Remark \ref{rem 2.2}, a field $K$ has Property $(*-)$ if and only if $K(\sqrt{-1})$ has Property $(*)$.

Having the results and conjectures discussed in this paper, a natural questions arises what can be said about the structure of fields with Property $(*)$ or $(*-)$. In particular, can one deduce from our results Conjectures \ref{con0}, \ref{con1}, \ref{con2}, \ref{con3} (at least in positive characteristic) or Podewski's conjecture?

Recall that a filed $K$ is bounded if its absolute Galois group is small, i.e., its absolute Galois group has only finitely many closed subgroups of any finite index; equivalently, for every natural number $n>0$, $K$ has only finitely many extensions of degree $n$ (up to isomorphism over $K$). A field is orderable if it can be equipped with some order making it an ordered field; equivalently, it is formally real (i.e., $-1$ is not a sum of squares).

Recall that a field $K$ is PAC (pseudo algebraically closed) if every absolutely irreducible variety over $K$ has a $K$-rational point. A field $K$ is PRC (pseudo real closed) if every absolutely irreducible variety over $K$ which has an $F$-rational point in every real closed field $F$ containing $K$ has a $K$-rational point. With such definitions each PAC field is PRC. We know that if $K$ is an orderable PRC field, then $K(\sqrt{-1})$ is (perfect) PAC \cite[Lemma A.1.1.3]{On}.

\begin{fact}\label{fact -1}
A PRC field is superrosy if and only if it is perfect and bounded.
\end{fact}
{\em Proof.} $(\leftarrow)$ was proved in \cite[Appendix A]{On} for orderable PRC fields and in \cite{Hr} for PAC fields. To see the converse, note that the same argument as in \cite[Theorem 5.6.5]{Wa} yields boundedness. Finally, suppose for a contradiction that $K$ is not perfect. Then, by Remark \ref{rem -11}, $[K^*:(K^*)^p]$ is infinite, which contradicts Fact \ref{fact rosy 2}. 
%guarantees that $K>K^p>K^{p^2}>\ldots$ is an infinite sequence of definable additive subgroups in which the quotients of successive terms are infinite. This contradicts superosiness, as $\uth(K)>\uth(K^p)>\dots$ would be an infinite decreasing sequence of ordinals.\hfill $\blacksquare$

%For the converse, suppose $K$ is a rosy PRC field which is not bounded. Then $K(\sqrt{-1})$ is still rosy (as a  field interpretable in $K$)

\begin{remark}\label{remark -2}
An orderable PRC field has the strict order property, so it is not simple. In particular, Conjecture \ref{con2} implies Conjecture \ref{con0}.
\end{remark}
{\em Proof.} Let $K$ be a formally real PRC field. We claim that any element which is a sum finitely many squares is a sum of two squares. For this, consider any non-zero element $a \in K$ which is a sum of finitely many squares. Then the polynomial $x^2+y^2-a$ has a zero in any real closed field containing $K$. Since $x^2+y^2-a$ is absolutely irreducible (e.g. by the Eisenstein criterion) and $K$ is PRC, this polynomial has a zero in $K$, and so $a$ is a sum of two squares in $K$.

Now, define the relation $\leq $ on $K$ by 
$$x \leq y \iff y-x\;\; \mbox{is a sum of squares of finitely many elements of}\;\; K.$$

We claim that $\leq$ is a partial order with an infinite chain. The fact that $\leq$ is antisymmetric follows from the fact that $K$ is formally real. Transitivity is clear. The existence of an infinite chain can be seen as follows. Take any $a \ne 0$. Then $a^2<2a^2<3a^2<\ldots$. 

Since, by the first paragraph of the proof, $x \leq y$ if and only if  $y-x$ is a sum of two squares, we see that $\leq$ is definable. \hfill $\blacksquare$

\begin{fact}\label{fact -3}
A non separably and non real closed PRC field does not have NIP.  In particular, Conjecture \ref{con2} implies Conjecture \ref{con1}.
\end{fact}
{\em Proof.} The main result of \cite{Du} says a non separably closed PAC field does not have NIP. Now, consider a non separably and non real closed PRC field. Suppose for a contradiction that it has NIP. Then $K(\sqrt{-1})$ still has NIP (as it is interpretable in $K$) and it is PAC but not separably closed, a contradiction. \hfill $\blacksquare$\\

%Since bounded orderable PRC fields and bounded perfect PAC fields are known to be superrosy, Corollay \ref{main corollary 1} implies that they have Property $(*-)$. In fact, in the case of PAC fields, $-1$ is a sum of squares, so the same property holds for all residue fields, and hence bounded perfect PAC fields have Property $(*)$. 

By \cite{Ef}, perfect PAC fields have Property $(*)$, so, by \cite[Lemma A.1.1.3]{On}, PRC fields have Property $(*-)$.
One could ask whether for infinite fields Property $(*-)$ [or $(*)$] implies that the field in question is PRC (in the non formally real case, PRC can be replaced by PAC). 
By virtue of our results and
%the fact that PAC non algebraically closed (so also PRC non real closed) fields do not have NIP (see \cite{Du}), 
Fact \ref{fact -3},
the positive answer would imply Conjectures \ref{con0}, \ref{con1}, \ref{con2} and \ref{con3} in positive characteristic. Unfortunately the answer is negative, which we briefly explain now.  In \cite{Ef} and \cite{GeJa}, the Hasse principle for Brauer groups is considered. It is shown in \cite[Theorem 3.4]{Ef} that if $K$ is a perfect PAC field, then any extension $F$ of $K$ of relative
transcendence degree 1 satisfies the Hasse principle for the Brauer groups. On the other hand, it is shown in \cite[Theorem 4.1]{Ef} that whenever K is a perfect field such that the Hasse principle for the
Brauer groups holds for all extensions $F$ of $K$ of relative transcendence degree 1, then the field $K$ has Property $(*)$. It was also asked in \cite[Question 4.2]{Ef} whether a non formally real infinite perfect field $K$ such that the Hasse
principle for the Brauer groups holds for all extensions $F$ of $K$ of relative
transcendence degree 1 is necessarily PAC? In \cite{GeJa}, a counter-example was constructed. The authors introduced the class of the so-called weakly PAC fields, they proved that weakly PAC fields are non formally real and that (assuming perfectness) they satisfy the Hasse principle for the
Brauer groups for all extensions of relative transcendence degree 1, and they constructed perfect weakly PAC fields which are not PAC. In particular, perfect weakly PAC fields have Property $(*)$. 

The above discussion leads to the following questions in our context.

\begin{question}
Does there exist a weakly PAC but not PAC superrosy [resp. supersimple] field?
\end{question}

The positive answer would refute Conjecture \ref{con2} [resp. Conjecture \ref{con0}]. The negative answer would support (but not prove) these conjectures.

\begin{question}\label{question 2}
Let $K$ be an infinite perfect field with NIP satisfying Property $(*)$ [or $(*-)$]. 
%Is it true that either there is a non-trivial definable valuation on $K$, or $K$ is algebraically or real closed?
Is it true that $K$ is either algebraically or real closed?
\end{question}

Applying the trick with adding $\sqrt{-1}$, we see that both versions of this question are equivalent.
By Corollaries \ref{main corollary 1} and \ref{main corollary 2}, the positive answer would imply Conjectures \ref{con1} and \ref{con3} in positive characteristic.

\begin{question}
Does there exist a non algebraically closed, perfect weakly PAC field which [is superrosy and] has NIP? 
\end{question}

The positive answer would yield the negative answer to Question \ref{question 2}. The positive answer to the extended version would refute Conjecture \ref{con1}. Notice that since by \cite{Du} we know that non separably closed PAC fields do not have NIP, if the answer to the above question was positive, the witness field would have to be weakly PAC but not PAC.

In any case, in order to find counter-examples to Conjectures \ref{con0},  \ref{con1}, \ref{con2} or \ref{con3}, one could try to construct  suitable 
%non formally real, non PAC fields satisfying Property $(*)$; 
non PRC fields with Property $(*-)$;
possible candidates could be among perfect weakly PAC but not PAC fields. Or, try to prove the conjectures by showing that there are no such fields.

A generalization of Corollary \ref{main corollary 1} to the characteristic zero case is an open problem. Of course, Thereom \ref{main theorem 1} yields the divisibility of the value groups, but the problem is with the residue fields.

\begin{main conjecture}\label{added conjecture}
Every non-trivial valuation on a superrosy field has divisible value group and either algebraically or real closed residue field.
\end{main conjecture}

Another idea of attacking Conjecture \ref{con1} by means of  Corollary \ref{main corollary 1} (or rather Conjecture \ref{added conjecture}), suggested by E. Hrushovski, is to try to show that a superrosy field with NIP is the residue field with respect to some non-trivial valuation on another superrosy field with NIP, and use Corollary \ref{main corollary 1} (or Conjecture \ref{added conjecture} in the zero characteristic case) to conclude that the original field is either algebraically or real closed.

As was mentioned in the introduction, K. Dupont has a different approach to Conjecture \ref{con3} (and so also to Conjecture \ref{con1}).  After adding $\sqrt{-1}$ to the field, the goal is to show that either for every finite extension $L$ of $K$ and for every prime number $n$, $(L^*)^n=L^*$ (as then $K$ is algebraically closed by Lemma \ref{lem 2.4}), or there is a non-trivial definable valuation on $K$.  Assume that the first possibility fails for some $L$ and $n$. Put $T=(L^*)^n$. It follows from the assumptions and Remark \ref{rem -11} that $L$ is perfect, so $n \ne \Char(L)$. The first goal is to show that ${\cal O}_T$ is a definable valuation ring on $L$. Further, K. Dupont deduced from \cite{Ko} that ${\cal O}_T$ is non-trivial if and only if  the family of sets $\{aT+a: a \in L^*\}$ is a subbasis of a $V$-topology on $L$. The second goal of her project is to show (using the NIP assumption) that the last condition holds.

\noindent
{\bf Address:}\\
Instytut Matematyczny, Uniwersytet Wroc\l awski,\\
pl. Grunwaldzki 2/4, 50-384 Wroc\l aw, Poland.\\[3mm]
{\bf E-mail address:}  kkrup@math.uni.wroc.pl

\end{document}